\documentclass[12pt]{amsart}
\usepackage{amscd,amssymb}
\usepackage[arrow,matrix]{xy}
\usepackage[plainpages,backref,urlcolor=blue]{hyperref}

\theoremstyle{plain}
\newtheorem{thm}[subsection]{Theorem}

\newtheorem{prop}[subsection]{Proposition}
\newtheorem{cor}[subsection]{Corollary}

\theoremstyle{definition}
\newtheorem{rk}[subsection]{Remark}
\newtheorem{definition}[subsection]{Definition}
\newtheorem{ex}[subsection]{Example}
\newtheorem{question}[subsection]{Question}

\numberwithin{equation}{section}
\newcommand{\RR}{{\mathcal R}}

\newcommand{\A}{{\mathcal A}}

\newcommand{\M}{{\mathcal M}}

\newcommand{\LL}{{\mathcal L}}

\newcommand{\V}{{\mathcal V}}

\newcommand{\al}{{\alpha}}
\newcommand{\be}{{\beta}}

\newcommand{\Z}{\mathbb{Z}}
\newcommand{\Q}{\mathbb{Q}}

\newcommand{\C}{\mathbb{C}}
\newcommand{\K}{\mathbb{K}}
\newcommand{\PP}{\mathbb{P}}
\newcommand{\FF}{\mathbb{F}}
\newcommand{\T}{\mathbb{T}}

\DeclareMathOperator{\Hom}{Hom}

\DeclareMathOperator{\coker}{coker}
\DeclareMathOperator{\ab}{ab}
\DeclareMathOperator{\id}{id}

\newcommand{\surj}{\twoheadrightarrow}

\begin{document}

\title [Finite covers, jump loci, formality, and multinets]
{Finite Galois covers, cohomology jump loci, formality properties, and multinets}

\author[Alexandru Dimca]{Alexandru Dimca$^1$}
\address{  Laboratoire J.A. Dieudonn\'e, UMR du CNRS 6621,
                 Universit\'e de Nice Sophia-Antipolis,
                 Parc Valrose,
                 06108 Nice Cedex 02,
                 France}
\email{dimca@unice.fr}

\author[Stefan Papadima]{Stefan Papadima$^2$}
\address{Institute of Mathematics Simion Stoilow, 
P.O. Box 1-764,
RO-014700 Bucharest, Romania}
\email{Stefan.Papadima@imar.ro}
\thanks{$^1$ Partially supported by the French-Romanian Programme LEA Math-Mode
and  ANR-08-BLAN-0317-02 (SEDIGA)} 
\thanks{$^2$ Partially supported by the French-Romanian Programme LEA Math-Mode
and CNCSIS ID-1189/2009-2011 (Romania).}



\begin{abstract}
We explore the relation between cohomology jump loci in a finite Galois cover,
formality properties and algebraic monodromy action. We show that the jump loci of the base and 
total space are essentially the same, provided the base space is $1$-formal and
the monodromy action in degree $1$ is trivial. We use reduced multinet structures on 
line arrangements to construct components of the first characteristic variety
of the Milnor fiber in degree $1$, and to prove
that the monodromy action is non-trivial in degree $1$. For an arbitrary line arrangement,
we prove that the triviality of the monodromy in degree $1$ can be detected in a precise way,
by resonance varieties.
\end{abstract}

\maketitle


\section{Introduction and statement of results} \label{sec:intro}

A homogeneous, degree $d$ polynomial $Q\in \C [x_0, \dots, x_n]$ defines a
hypersurface in $\PP^n$, $V=V(Q)= \{ Q=0\}$, its complement $M=M(Q)$, the
Milnor fiber in $\C^{n+1}$, $F=F(Q)= \{ Q=1\}$, and a $d$-fold cyclic Galois cover,
$p:F\to M$, with geometric monodromy $h:F\to F$. In spite of the abundance of available
algebro-geometric and topological methods, the analysis of the algebraic monodromy
action (over $\C$), $h^*:H^*(F)\to H^*(F)$, remains a challenge, as documented by an
extensive literature. Even when $Q$ completely decomposes into distinct linear factors, 
the additional combinatorial tools from the theory of hyperplane arrangements have not succeeded
to elucidate the subject. 

Given a connected, finite CW-complex $X$, one may consider its associated {\em characteristic
varieties} $\V^i_m(X)$ (alias its Green-Lazarsfeld sets, when $X$ is a smooth projective variety,
see \cite{GL} for the original setting involving line bundles with trivial Chern class, and \cite{Beau} for
the reformulation in terms of local systems).
They are defined as jump varieties in degree $i$, for the cohomology of $X$ with coefficients in
rank $1$ local systems, and live in $\T (X):= \Hom (\pi_1 (X), \C^*)$. Their `infinitesimal'
approximations, the {\em resonance varieties} $\RR^i_m (X)$ associated to the cohomology ring of $X$,
are homogeneous subvarieties of $H^1(X, \C)$. See Section \ref{sec:two} for the precise definitions.
Encouraged by recent progress on jump loci, we aim in this note to start their study for Milnor
fibers, and relate them to the algebraic monodromy.

\subsection{} \label{ss11}

For a `general' hyperplane arrangement in $\PP^n$, it is well-known that $h^* =\id$ in
degree $*<n$. See e.g. \cite{DP}, \cite{CS1}, \cite[Chapter 6]{D2}, \cite{L1}, \cite{MP},
for various (classes of) examples of this kind. Guided by these examples, we begin our study
in Section \ref{sec:two} with jump loci in arbitrary covers, $p:F\to M$, with finite Galois group
$C$, with emphasis on  implications of the (partial) triviality of the algebraic monodromy
action of $C$ on $H^*(F)$. 

We start by showing that $p^*: \Hom (\pi_1 (M), \C^*)\to \Hom (\pi_1 (F), \C^*)$ sends $\V^i_m(M)$
into $\V^i_m(F)$, for all $i$ and $m$. For the rest of Section \ref{sec:two}, we fix $q\ge 1$ and
suppose that $C$ acts trivially on $H^*(F)$, for $*\le q$. We first note that under this assumption
the isomorphism $p^*: H^1(M, \C) \stackrel{\sim}{\longrightarrow}  H^1(F, \C)$ identifies 
$\RR^i_m (M)$ with $\RR^i_m (F)$, for all $i\le q$ and all $m$. Unfortunately, the proof 
does not work for characteristic varieties.

This difficulty may be bypassed by resorting to a recent result from \cite{DPS08}.
This result (recorded here as Theorem \ref{thm:expf}) establishes a natural bijection between
the {\em non-translated} irreducible components of $\V^1_m(X)$ (i.e., those containing the
trivial local system $1$) and the 
irreducible components of $\RR^1_m(X)$, for all $m$, under a so-called {\em $1$-formality}
hypothesis on the connected CW-complex $X$.

We are thus naturally led to examine partial formality properties in finite Galois covers.
Their definition is inspired from D.~Sullivan's \cite{Sul} homotopy theory for commutative
differential graded algebras in characteristic zero. To any such DGA, $(A^*, d_A)$, one may
associate its cohomology algebra $H^*(A)$, viewed as a DGA with trivial differential. For a
fixed $q\ge 1$, $(A^*, d_A)$ is called {\em $q$-formal} if it can be connected to
$(H^*(A), d=0)$ by a zigzag of DGA maps, each inducing in cohomology an isomorphism in degree
up to $q$, and a monomorphism in degree $q+1$. A connected polyhedron $X$ is called 
$q$-formal if Sullivan's DGA of PL forms on $X$ is $q$-formal. See Definition \ref{def:partial}
and Remark \ref{rem:simple} for another, equivalent, approach to this notion. For $q=\infty$,
one recovers Sullivan's celebrated notion of {\em formality}, which implies $q$-formality, for
all $q$. The $1$-formality of $X$ depends only on $\pi_1(X)$. See Section \ref{sec:two} for more
details.

Interesting examples abound. Hyperplane arrangement complements are formal \cite{Br}.
If a smooth quasi-projective variety $X$ has vanishing Deligne weight filtration $W_1H^1(X)$
(which happens e.g. for a projective hypersurface complement), then $X$ is $1$-formal \cite{M}.
In particular, pure braid groups are $1$-formal. This was extended to pure welded braid groups
in \cite{BP}. Finitely generated Artin groups are also $1$-formal \cite{KM}. Under certain
trivial monodromy assumptions, Artin kernels are $1$-formal, too \cite{PS-t}.

The interplay between $1$-formality and algebraic monodromy in fibrations was recently
investigated in \cite{PS-m}. We pursue this approach for finite Galois covers with trivial 
algebraic monodromy up to degree $q$,  proving that $M$ is $q$-formal if and only if $F$ is
$q$-formal. Here is our main result from Section \ref{sec:two}.

\begin{thm}[Theorem \ref{thm:monoch}]
\label{thma}
Let $p:F\to M$ be a finite Galois cover with group $C$. Assume that $C$ acts trivially on $H^1(F)$,
and $M$ is $1$-formal. Then $p^*: \T (M)\to \T (F)$ gives a natural bijection between 
non-translated components in $\V^1_m(M)$ and $\V^1_m(F)$, for all $m$.
\end{thm}

\subsection{} \label{ss12}

In Section \ref{sec:3}, we specialize our discussion to the Galois cover associated to
a line arrangement in $\PP^2$, $\A =\{ L_1, \dots, L_d \}$, with $L_i$ given by the equation
$f_i=0$. It follows from Theorem \ref{thm:expf} that the non-translated components of
$\V^1_m(M)$ are in bijection with the components of $\RR^1_m(M)$, for all $m$.
At the same time, it is well-known that $\RR^1_m(M)$ depends only on the combinatorics of the
arrangement, as follows from basic work by Orlik and Solomon \cite{OS}. A nice recent result 
of Falk and Yuzvinsky \cite{FY} gives a combinatorial parametrization for {\em global} components
of $\RR^1_1(M)$ (i.e., components not coming from a proper subarrangement).

The key combinatorial notion from \cite{FY} is that of a {\em multinet} supported by $\A$.
It involves two data: a multiplicity function $\mu$ on the set $\{ 1, \dots, d\}$ 
with values in $\Z_{>0}$, and a partition of this set, $\A_1 \cup \dots \cup \A_k$, 
with at least $3$ elements, satisfying certain axioms. The first axiom requires that
$\sum_{i\in \A_j} \mu (i)=e$, independently of $j$. Starting with a multinet structure on $\A$,
the authors of \cite{FY} construct an {\em admissible map} (in the sense of Arapura \cite{A}),
$f:M\to S$, where $S=\PP^1 \setminus \{ k \; {\rm points}\}$, and the global component
$f^*(H^1(S))$ of $\RR^1_1(M)$, whose corresponding non-translated component of $\V^1_1(M)$
is $W:= f^*(\T(S))$. They also prove that every global component of $\RR^1_1(M)$ arises in this way.

We will say that a multinet is {\em reduced} if $\mu$ takes only the value $1$. Classical nets 
provide reduced examples. There are also situations when $\A$ supports no multinet, simply because
there is no global component of $\RR^1_1(M)$; see Example 10.5 from \cite{S2}. In a recent preprint
\cite{Y}, Yuzvinsky shows that $k$ is $3$ or $4$, for any multinet. The only known example with $k=4$
comes from the Hesse pencil.

Recall from Theorem \ref{thma} that $p^*(W)$ is a non-translated component of $\V^1_1(F)$, if
$h^*=\id$ on $H^1(F)$, where $W$ is the  non-translated component of $\V^1_1(M)$ associated
to an arbitrary multinet on $\A$. Our second main result in this note is derived in
Section \ref{sec:3}. It establishes a connection between the existence of reduced multinet
structures on a line arrangement, and the triviality of the algebraic monodromy in degree $1$.

\begin{thm}
\label{thmb}
Let $W$ be the  non-translated component of $\V^1_1(M)$ associated
to a reduced multinet on $\A$. Then there is a non-translated component of $\V^1_1(F)$, $W'$,
strictly containing $p^*(W)$. In particular, $h^*\ne \id$, on $H^1(F)$.
\end{thm}

See Theorem \ref{thm3} and Corollary \ref{cor:charnet} for more precise statements. 
We point out that the existence of a non-reduced multinet is not enough to infer the
non-triviality of the monodromy; see Remark \ref{rk2}$(ii)$.

In Section \ref{sec:mhs}, we examine Galois covers associated to plane projective curves,
and exploit the mixed Hodge structure. When $p:F\to M$ comes from 
an arbitrary line arrangement, we show that $h^*=\id$ on $H^1(F)$ if and only if 
$p^*(\RR^1_1(M))=\RR^1_1(F)$, thereby improving Theorem \ref{thma} in this case.
The proof relies on a general result that relates the mixed Hodge structure on $H^*(F)$ and 
the algebraic monodromy action; see Theorem \ref{thm2}. Finally, we extend the natural 
bijection from Theorem \ref{thma}, for arbitrary plane curves. In this situation, we show that,
when $h^*=\id$ on $H^1(F)$, $p^*$ essentially identifies all components 
in $\V^1_1(M)$ and $\V^1_1(F)$, of dimension at least $2$. See Proposition \ref{prop1}
for a more precise statement.

\bigskip

In the recent preprint \cite{BDS}, the relation between nets and monodromy of the Milnor fiber of a line arrangement is considered from an alternative viewpoint, see especially \cite[Theorem 1]{BDS}.

\section{Finite Galois covers} \label{sec:two}

In this section, we begin our analysis of cohomology jump loci in a purely topological context.
Let $X$ be a connected CW-complex with finite $1$-skeleton, and (finitely generated)
fundamental group $G=\pi_1(X)$. For a characteristic zero field $\K$, denote by
$\T(X, \K)=\Hom (G, \K^*)= \Hom (G_{\ab}, \K^*)$ the {\em character torus}.
(When there is no explicit mention about coefficients, we mean that $\K=\C$.) Clearly,
the algebraic group $\T(X, \K)$ is the direct product with a finite group of the
{\em connected character torus}, $\T^0(X, \K)=\Hom (G_{\ab}/ {\rm torsion} , \K^*)= (\K^*)^{b_1(G)}$.

The characteristic varieties
\begin{equation} 
\label{eq1}
\V^i_m(X, \K)=\{\LL \in \T(X, \K) \mid \dim_{\K} H^i(X, \LL)\ge m\}
\end{equation}
are defined for $i\ge 0$ and $m>0$. If $X$ has finite $q$-skeleton ($q\ge 1$),
they are Zariski closed subsets of the character torus, for $i\le q$ and $m>0$.
Their approximations, the resonance varieties
\begin{equation} \label{e4}
\RR^i_m(X, \K)=\{\alpha \in H^1(X,\K) \mid \dim_{\K} H^i(H^{\bullet}(X,\K), \mu_{\alpha})\ge m\}
\end{equation}
(where $\mu_{\alpha}$ denotes right-multiplication by $\alpha$ in the cohomology ring of $X$)
are defined for $i\ge 0$ and $m>0$. If $X$ has finite $q$-skeleton ($q\ge 1$),
they are homogeneous, Zariski closed subsets of the affine space $H^1(X, \K)= \K^{b_1(G)}$.
We mean that $i=1$, whenever the upper index $i$ is missing. In this case, both $\V_m(X, \K)$ and
$\RR_m(X, \K)$ depend only on $G$, for all $m$.

Let $p:F\to M$ be a connected cover of a complex $M$ having finite $1$-skeleton, 
with  Galois group $C$.

\begin{prop}
\label{prop:galjump}
If $C$ is finite, then 
\[
p^*(\V^i_m(M, \K))\subseteq \V^i_m(F, \K) \quad {\rm and}\quad 
p^*(\RR^i_m(M, \K))\subseteq \RR^i_m(F, \K)\, ,
\]
for all $i\ge 0$ and $m>0$.
\end{prop}

\begin{proof}
Let $\LL$ be an arbitrary right $\pi_1(M)$-module. In the Hochschild-Serre spectral
sequence of an arbitrary connected Galois cover of a CW-complex, $p:F\to M$,
\[
E^2_{st}= H_s(C, H_t(F, p^* \LL)) \Rightarrow H_{s+t}(M, \LL)\, ,
\]
see e.g. Brown's book \cite{B}, Proposition 5.6 on p. 170 and Theorem 6.3 on p. 171.
When $C$ is finite and $\LL$ is a $\K$-representation in characteristic zero, $E^2_{st}=0$,
for all $s>0$ and $t\ge 0$, and the spectral sequence collapses to the isomorphism
$H_{\bullet}(M, \LL)= H_{\bullet}(F, p^*\LL)_C$, where $(\cdot)_C$ denotes coinvariants.
By duality, we obtain an injection, 
$p^{\bullet}: H^{\bullet}(M, \LL)\hookrightarrow H^{\bullet}(F, p^*\LL)$, for any $\LL \in \T(M, \K)$. 
Our claim on characteristic varieties follows.

In the particular case when the character (local system) $\LL$ is trivial, 
we may identify the induced algebra map in cohomology, 
$p^{\bullet}: H^{\bullet}(M, \K)\rightarrow H^{\bullet}(F, \K)$,
with the inclusion of fixed points,
\begin{equation}
\label{eq:fix}
p^{\bullet}: H^{\bullet}(F, \K)^C \hookrightarrow H^{\bullet}(F, \K) \,.
\end{equation}
For $\alpha \in H^{1}(F, \K)^C$, the monodromy $C$-action on $H^{\bullet}(F, \K)$
clearly gives rise to a $C$-action on the chain complex $(H^{\bullet}(F, \K), \mu_{\alpha})$,
with fixed subcomplex $(H^{\bullet}(M, \K), \mu_{\alpha})$. Since $C$ is finite and char $\K =0$,
we obtain an inclusion, 
$H^*(H^{\bullet}(M, \K), \mu_{\alpha}) \hookrightarrow H^*(H^{\bullet}(F, \K), \mu_{\alpha})$.
This proves our claim on resonance varieties.
\end{proof}

We continue by examining the algebraic monodromy action, and the implications of its
triviality on cohomology jump loci. We keep the hypotheses from Proposition \ref{prop:galjump}.

\begin{cor}
\label{cor:monores}
If $C$ acts trivially on $H^i(F, \K)$, for $i\le q$, where $q\ge 1$, then
$p^{*}: H^{*}(M, \K)\rightarrow H^{*}(F, \K)$ is an isomorphism for $*\le q$ and a 
monomorphism for $*=q+1$. In degree $*=1$, $p^*$ identifies $\RR^i_m(M, \K)$ with
$\RR^i_m(F, \K)$, for all $i\le q$ and $m>0$.
\end{cor}

\begin{proof}
Everything follows from identification \eqref{eq:fix}.
\end{proof}

To obtain a similar result for characteristic varieties, we need to discuss first
partial formality (in the sense from \cite{Mac}). This notion can be naturally extracted from
Sullivan's foundational work in rational homotopy theory \cite{Sul} (see also Morgan \cite{M}).
Sullivan associates to a space $X$, in a natural way, a $\K$-DGA, $\Omega^*_{PL}(X, \K)$; see
also Bousfield and Gugenheim \cite{BG}. This is an analog of the usual de Rham DGA of a manifold.
He proves that $H^{\bullet}\Omega^*_{PL}(X, \K)=H^{\bullet}(X, \K)$, as graded $\K$-algebras.

Sullivan goes on defining, for each $1\le q\le \infty$, an important class of $\K$-DGA's: the
{\em $q$-minimal} ones (simply called {\em minimal}, for $q=\infty$). A DGA map,
$f: (A^*, d_A)\to (B^*, d_B)$, is  a {\em $q$-equivalence} 
(or, {\em quasi-isomorphism}, for $q=\infty$) if the induced map in cohomology is an isomorphism
up to degree $q$, and an injection in degree $q+1$. Given any homologically connected DGA 
$(A^*, d_A)$ (i.e., such that $H^0(A, d_A)=\K$), he shows that there is a $q$-minimal DGA,
$(\M_q, d)$, together with a $q$-equivalence, $\varphi :(\M_q, d)\to (A, d_A)$. Moreover,
these two conditions uniquely determine $(\M_q, d)$, up to DGA isomorphism. The DGA
$(\M_q, d)$ is called the {\em $q$-minimal model} of $(A, d_A)$ (simply the {\em minimal model},
for $q=\infty$), and will be denoted by $\M_q(A, d_A)$ (respectively by $\M(A, d_A)$, for $q=\infty$).
For $1\le r\le q\le \infty$, $\M_r(A, d_A)$ is the sub-DGA of $\M_q(A, d_A)$ generated by the elements
of degree at most $r$. 

Sullivan calls a homologically connected $\K$-DGA $(A, d_A)$ formal if 
$\M(A^*, d_A)= \M(H^*(A, d_A), d=0)$ , up to DGA isomorphism. A path-connected space $X$
is formal (over $\K$) if the DGA  $\Omega^*_{PL}(X, \K)$ is formal. Since the appearance
of \cite{DGMS}, where the authors prove the formality of compact Kahler manifolds, this notion
played a key role in rational homotopy theory and its applications.
The following notion of partial formality was investigated by Macinic in \cite{Mac}.

\begin{definition}
\label{def:partial}
A $\K$-DGA $(A^*, d_A)$ with $H^0(A, d_A)=\K$ is $q$-formal ($1\le q\le \infty$)
if $\M_q (A^*, d_A)= \M_q (H^*(A, d_A), d=0)$, up to DGA isomorphism. A path-connected space $X$
is $q$-formal (over $\K$) if the DGA $\Omega^*_{PL} (X, \K)$ is $q$-formal.
\end{definition}

\begin{rk}
\label{rem:simple}
A couple of useful simple properties follow directly from Definition \ref{def:partial},
via the discussion preceding it. $X$ is $\infty$-formal if and only if $X$ is formal in the sense 
of Sullivan. For $1\le r\le q\le \infty$, $q$-formality implies $r$-formality. If
$p:F\to M$ is a continuous map between path-connected spaces, inducing in $\K$-cohomology an
isomorphism up de degre $q$, and a monomorphism in degree $q+1$, then 
$F$ is $q$-formal if and only if $M$ is $q$-formal, over $\K$. Taking $M$ to be a point, it
follows that $F$ is $q$-formal, if $b_i(F)=0$, for $i\le q$. See also the survey \cite{PSs} for more details on partial formality.
\end{rk}

\begin{cor}
\label{cor:galf}
Let $p:F\to M$ be a finite Galois cover, as in Proposition \ref{prop:galjump}. Assume that $C$
acts trivially on $H^i(F, \K)$, for $i\le q$ ($q\ge 1$). Then 
$F$ is $q$-formal if and only if $M$ is $q$-formal, over $\K$.
\end{cor}

\begin{proof}
Combine Corollary \ref{cor:monores} and Remark \ref{rem:simple}.
\end{proof}

\begin{cor}
\label{cor:galarr}
Let $p:F\to M$ be the cyclic Galois cover associated to a hyperplane arrangement 
in $\PP^n$ ($n\ge 2$). If $h^*: H^*(F)\to H^*(F)$ is the identity, for $*<n$, then
the Milnor fiber $F$ is a formal space.
\end{cor}

\begin{proof}
By Corollary \ref{cor:galf}, $F$ is $(n-1)$-formal. Up to homotopy, $F$ is a 
CW-complex of dimension at most $n$. These two facts together imply the formality of $F$; 
see \cite{Mac}.
\end{proof}

Let $p:X \to K(G,1)$ be the classifying map of a connected complex with fundamental group
$G=\pi_1(X)$. By Remark \ref{rem:simple}, $X$ is $1$-formal if and only if $K(G, 1)$ is $1$-formal.
When this happens, we simply say that the group $G$ is $1$-formal. It is implicit in \cite{Sul}
that a finitely generated group $G$ is $1$-formal if and only if its Malcev Lie algebra
(in the sense of Quillen \cite{Q}) is a quadratic complete Lie algebra.

The $1$-formality property has many interesting consequences. At the level of cohomology
jump loci, these may be described as follows. Let $X$ be a connected CW-complex with finite
$1$-skeleton. For each $m>0$, denote by $\breve{\V}_m(X) \subseteq \T^0(X)$ the union of all
non-translated irreducible components of $\V_m(X)$. Let $\exp \colon H^1(X) \surj \T^0(X)$ be
the exponential map.

\begin{thm}[\cite{DPS08}] 
\label{thm:expf}
If $X$ is $1$-formal, then $\RR_m(X)$ is a finite union of linear subspaces of $H^1(X)$,
defined over $\Q$, and $\breve{\V}_m(X) = \exp (\RR_m(X))$, for all $m>0$. More precisely,
if $\{ U_i\}$ are the irreducible components of $\RR_m(X)$, then $\{ \exp (U_i)\}$ are
the irreducible components of $\breve{\V}_m(X)$.
\end{thm}

\begin{thm}
\label{thm:monoch}
Let $p:F\to M$ be a connected, finite Galois cover with group $C$, where $M$ is
a CW-complex with finite $1$-skeleton. Assume that $C$ acts trivially on $H^1(F)$,
and $M$ is $1$-formal. Then the following hold. 

The morphism $p^*: \T^0(M)\surj \T^0(F)$ is a surjection with finite kernel.
The space $F$ is $1$-formal. The irreducible components of $\breve{\V}_m(F)$
are $\{ p^*(W)\}$, where $W$ runs through the irreducible components of $\breve{\V}_m(M)$,
for all $m>0$.
\end{thm}

\begin{proof}
The first assertion follows from the fact that $p^*: H^1(M)\to H^1(F)$ is an isomorphism.
The second claim is implied by Corollary \ref{cor:monores}, via Remark \ref{rem:simple}.
Again by Corollary \ref{cor:monores}, the isomorphism $p^*: H^1(M)\to H^1(F)$ identifies the
irreducible components of $\RR_m(F)$ with those of $\RR_m(M)$. Let 
$\{ U_i\}$ be the irreducible components of $\RR_m(M)$. We infer from Theorem \ref{thm:expf} 
that $\{ \exp (U_i)\}$ are
the irreducible components of $\breve{\V}_m(M)$, and 
$\{ \exp \circ p^*(U_i)\}= \{ p^* \circ\exp (U_i)\}$ are those of $\breve{\V}_m(F)$, 
which proves the last claim.
\end{proof}

We point out that the triviality assumption on the algebraic monodromy from 
Corollary \ref{cor:galf} is necessary, since in general the formality properties 
of the base space are not inherited by the total space of a finite Galois cover.
The following simple example arose from conversations with Alex Suciu.

\begin{ex}
\label{ex:bigu}
Let $\FF$ be the free group generated by $x$ and $y$. Let $K$ be the quotient of $\FF$ by
$3$-fold commutators. Let $\varphi$ be the order $2$ automorphism of $K$ induced by
$x\mapsto x^{-1}$, $y\mapsto y^{-1}$. Using $\varphi$, we may construct the semidirect product
$G=K\rtimes (\Z/2 \Z)$. Let $M$ be the connected, finite $2$-complex associate to 
a presentation of $G$. Consider the $2$-cover, $p:F\to M$, with $\pi_1(F)=K$. Clearly,
$h^*= -\id$ on $H^1(F)$. Consequently, $b_1(M)=0$. By Remark \ref{rem:simple}, $M$ is 
$1$-formal. Hence, $M$ is formal, see \cite{Mac}. On the other hand, $K$ (hence $F$) is not
$1$-formal, as noticed by Morgan in \cite{M}.
\end{ex}

\begin{question} \label{q}
Is the Milnor fiber of an arbitrary hyperplane arrangement $1$-formal?
(Note that smooth affine varieties need not be $1$-formal; see \cite[Proposition 7.2]{DPSqk}.)
\end{question}

\section{`Exceptional' arrangements} \label{sec:3}

Let $\A =\{ L_1, \dots, L_d\}$ be a line arrangement in $\PP^2$,
with associated $d$-fold cyclic Galois cover $p:F\to M$.
Let $\A_1 \cup...\cup \A_k$ be a partition of the set $\{1,2,\dots,d\}$ into $k\geq 3$ subsets
of the same cardinality $e>0$.
Let $f_i=0$ be an equation for the line $L_i$ for $i=1,...,d$, and set $Q_j=\prod_{i\in \A_j}f_i$,
for $j=1,...,k$. Clearly, $Q=Q_1 \cdots Q_k$ is the defining polynomial of $\A$.

\begin{thm} \label{thm3}
With the above notation, assume that the vector space $<Q_1,...,Q_k>$ of degree $e$ 
homogeneous polynomials has dimension $2$. 
Then the following hold.

\medskip

\noindent (i) $h^*:H^1(F) \to H^1(F)$ is not the identity; more precisely, the eigenspace 
$H^{1}(F)_{\lambda}$ has dimension at least $k-2$, for any $\lambda$ with $\lambda^k=1$.

\medskip

\noindent (ii) Assume in addition that the pencil $f:M \to S$ given by $f=(Q_1,Q_2)$, 
(where $S$ is obtained from $\PP^1$ by deleting $k$ points) has a connected generic fiber.
Let $I=f^*(H^1(S))$ be the corresponding maximal isotropic subspace in $H^1(M)$ 
(with respect to the cup product). Then $\dim I=k-1$ and
there is an admissible morphism (in the sense of \cite{A}), $f':F \to S'$, where
$\chi (S')<0$, such that
the corresponding subspace $J=f^{'*}(H^1(S'))$  in $H^1(F)$ has dimension 
at least $(k-1)^2$ and satisfies
$J \cap p^*(H^{1}(M))=  p^*(I) $.
\end{thm}

\begin{proof}
To prove $(i)$, consider the pencil $f:M \to S$ given by $f=(Q_1,Q_2)$. 
Let $B=\{b_1,...,b_k\}$ be 
the finite set such that $S=\PP^1\setminus B$. By the Stein factorization, 
there is a finite map $p_0:S_0 \to S$ and
a morphism $f_0:M \to S_0$ such that the generic fiber of $f_0$ is connected, 
and $f=p_0 \circ f_0 $.
Note that $S_0$ is a non-compact curve, and $p_0^*:H^1(S) \to H^1(S_0)$ 
is injective by Lemma 6.10 in \cite{DPS08}.
It follows that
$$\chi(S_0) \leq \chi(S) =2-k \leq -1.$$
By Arapura's work, see \cite{A}, \cite{D3}, it follows that $W=f_0^*(\T(S_0))$ is 
an irreducible component of
$\V_1(M)$ and that for any $\LL \in W$ one has
$$\dim H^1(M,\LL) \geq \dim W-1=b_1(S_0)-1 \geq b_1(S)-1=k-2.$$
Note also that $f^*(\T(S)) \subset W$.

On the other hand, a local system $\LL' \in \T(S)$ is determined by a family of complex numbers
$$(\lambda_1,...,\lambda_k) \in (\C^*)^k$$
satisfying $\lambda_1 \cdot...\cdot \lambda_k=1$, where $\lambda_j$ is the monodromy of the local
system $\LL'$ about the point $b_j$. For any $\lambda \in \C$ satisfying $\lambda^k=1$, 
we denote by $\LL'_{\lambda}$ the local system in $\T(S)$ corresponding to the choice
$$\lambda_1=...=\lambda_k=\lambda.$$
The pull-back local system $\LL_{\lambda}=f^*\LL'_{\lambda}$ is the unique local system on $M$
whose monodromy about each line is ${\lambda}$. By \cite{CS1}, \cite{D2} we know that
$$\dim H^{1}(F)_{\lambda}=\dim H^{1}(M, \LL_{\lambda}).$$
This dimension is at least $k-2$ since $\LL_{\lambda} \in W$, which proves claim $(i)$.

To prove  claim $(ii)$, we proceed as follows.
For $j=3,...,k$ there are unique complex numbers $\al_j,\be_j$ such that $Q_j=\al_jQ_1+\be_jQ_2$.
Consider the homogeneous polynomial $G(u,v)=uv\prod_{j=3,k}(\al_ju+\be_jv)$. 
Let $S=M(G)$ and $H=F(G)$
be the corresponding complement in $\PP^1$ and Milnor fiber in $\C^2$.

The map $g:F \to H$, $g(x,y,z)=(Q_1(x,y,z),Q_2(x,y,z))$ is well-defined 
(since $Q=Q_1\cdot ...\cdot Q_k$) and surjective. 
If the generic fiber of $g$ is connected, then we take $S'=H$ and $f'=g$. Otherwise, there is 
a Stein factorization $g=p' \circ f'$, where $p':S' \to H$ is a finite morphism and $f':F \to S'$
has connected generic fibers. Note that both $p'$ and $f'$ induce monomorphisms in cohomology, 
in degree $1$, by \cite[Lemma 6.10]{DPS08}. Using the same argument as in Part $(i)$, we may infer that
$\chi(S') \le \chi(H)<0$.

To show that $J \cap p^*(H^{1}(M))=  p^*(I)$, use the obvious inclusion
$J \cap p^*(H^{1}(M))\supseteq p^*(I) $, and the fact that $I$ is a maximal isotropic subspace 
in $H^1(M)$, see \cite{D4}. Finally, 
$\dim J = \dim H^1(S')\geq \dim H^1(H)=(k-1)^2$. This completes the proof.
\end{proof}

\begin{cor}
\label{cor:charnet}
Under the assumptions from Theorem \ref{thm3}$(ii)$, there exist components, 
$W\subseteq \V_1(M)$ and $W' \subseteq \V_1(F)$, both non-translated, 
with $\dim W'> \dim W$, whose tangent spaces at $1$
satisfy $p^*(T_1W)= T_1W' \cap p^*(H^1(M))$. In particular, the monodromy action on $H^1(F)$ is
non-trivial.
\end{cor}

\begin{proof}
Consider the admissible maps, $f:M\to S$ and $f':F\to S'$, from Theorem \ref{thm3}$(ii)$. Since
$\chi(S), \chi(S')<0$, we obtain non-translated components, $W=f^*(\T(S))$ and $W'=f'^*(\T(S'))$,
by Arapura theory \cite{A}. Clearly, $T_1W=I$ and $T_1W'=J$, using 
the notation from Theorem \ref{thm3}. All claims except the last (on monodromy) follow from
Theorem \ref{thm3}$(ii)$. The triviality of $h^*:H^1(F)\to H^1(F)$ would imply that $p^*(W)$
is an irreducible component of $\V_1(F)$, due to Theorem \ref{thm:monoch}. This contradicts
the strict inclusion $p^*(W)\subset W'$.
\end{proof}

\begin{cor}
\label{cor:subtle}
The admissible map
$f':F \to S'$ from Theorem \ref{thm3}$(ii)$  is not a rational pencil, i.e., 
$S'$ is not an open subset of $\PP^1$.
\end{cor}

\begin{proof} 
Consider the induced map on compactifications, $\widehat{p'} :\widehat{S'} \to \widehat{H}$,
to infer that $b_1(\widehat{S'})\ge b_1(\widehat{H})= (k-1)(k-2)>0$.
\end{proof}

\begin{rk} \label{rk2}

$(i)$ At this point, it seems worth recalling a couple of relevant facts from \cite{FY}.
Firstly, admissible maps may be constructed from arbitrary multinet structures, 
in the following way. Set $Q_j= \prod_{i\in \A_j} f_i^{\mu(i)}$, for $j=1, \dots, k$. 
Then the subspace $<Q_1, \dots, Q_k>$ is $2$-dimensional, and the associated map, 
$f=(Q_1, Q_2): M\to S=\PP^1 \setminus \{ k\; {\rm points}\}$, has connected generic fiber,
as needed in Theorem \ref{thm3}. Secondly, a degree $e$ pencil with connected fibers
(as defined in \cite{FY}), having $k\ge 3$ completely reducible fibers,
$Q_j=\prod f_i^{\mu(i)}$ (where all factors $f_i$ have degree $1$), for $j=1, \dots, k$, 
gives rise to a multinet structure with multiplicity function $\mu$, on the arrangement
consisting of the lines $\{ f_i=0\}$. Finally, every multinet arises in this way.
Clearly, the reduced multinets correspond to the case when all special fibers are reduced.
When such a reduced structure exists, $h^* \ne \id$, on $H^1(F)$, by Corollary \ref{cor:charnet}.
The existence of a reduced multinet also implies that the characteristic varieties of 
the Milnor fiber $F$ are more subtle
than those of a line arrangement complement, which are known to come from rational pencils; see
Corollary \ref{cor:subtle}.

\medskip

$(ii)$ If only non-reduced multinet structures exist, the monodromy action on $H^1(F)$ 
may well be trivial. Consider the $B_3$-arrangement, whose defining equation is given by
$Q= xyz (x^2-y^2)(y^2-z^2)(z^2-x^2)$. The pencil
$$<Q_1=x^2(y^2-z^2), Q_2=y^2(z^2-x^2), Q_3=z^2(x^2-y^2)>$$
gives rise to a multinet on $B_3$, with multiplicities equal to $2$ on $x$, $y$ and $z$,
and equal to $1$, otherwise. See \cite[Example 3.6]{FY}.
As shown in \cite{CS1}, the monodromy action on $H^1(F)$ is trivial.

\medskip
 
$(iii)$ We point out that the non-triviality of the monodromy action on $H^1(F)$ may  
be deduced from the existence of a non-reduced multinet structure, provided its multiplicity
function $\mu$ enjoys some special properties. Let us examine  the family of line arrangements
$\A_r$ from \cite[Example 4.6]{FY}, having defining polynomial
$$xyz(x^r-y^r)(y^r-z^r)(z^r-x^r)\,,\quad {\rm and}\quad d=3(r+1).$$
The corresponding pencil is 
$$<Q_1=x^r(y^r-z^r), Q_2=y^r(z^r-x^r), Q_3=z^r(x^r-y^r)>,$$
and the multiplicity function takes the value $r$ on on $x$, $y$ and $z$,
and $1$, otherwise. Assume that $r\equiv 1$ (modulo $3$). 
Let $\lambda$ be a primitive root of unity of order $3$. Then the local system 
$\LL_{\lambda}$ constructed
as in the proof of Theorem \ref{thm3}$(i)$ has monodromy  $\lambda$ about every line of $\A_r$
(since  $\lambda^r= \lambda$). Consequently, the same argument as in Theorem \ref{thm3}$(i)$
implies that the eigenspace $H^{1}(F)_{\lambda}$ has dimension at least $1$.

\medskip
 
$(iv)$ The inequality $\dim H^1(F)_{\lambda} \geq k-2$ from Theorem \ref{thm3} $(i)$ above is shown to be an equality under some additional condition in \cite[Theorem 1 (ii)]{BDS}.
 
\end{rk}

\begin{ex} \label{ex1}
For the $A_3$-arrangement and the Pappus configuration $(9_3)_1$, it follows from  \cite{CS1} 
that the inequality in Theorem \ref{thm3} $(i)$ is in fact an equality. This is the case whenever 
the generic fiber of
$g:F \to H$ is connected, i.e. $f'=g$, and the local system $\LL_{\lambda}$ is 1-admissible, 
which is equivalent to $\dim H^1(M, \LL_{\lambda})=\dim W-1$, see \cite{D5} for this equivalence.
\end{ex}

\begin{ex} \label{ex2}
Consider the Hesse arrangement  consisting of the 12 lines that occur in the 4 special fibers of  
the Hesse pencil $f=(x^3+y^3+z^3,xyz)$; see \cite{FY}. Here $d=12$, 
and there is a partition with $k=4$.
Theorem \ref{thm3} implies that in this case $H^{1}(F)_{-1}$, $H^{1}(F)_{i}$ and 
$H^{1}(F)_{-i}$   have each dimension at least $2$. Using \cite[Theorem 1 (ii)]{BDS}, one finds that all these dimensions are equal to $2$.
\end{ex}

\section{Monodromy action and mixed Hodge structure} \label{sec:mhs}
 
We begin this section with a general result on Milnor fibers, valid for a class of polynomials
that includes the case of line arrangements.
Let $Q \in \C[x_0,...,x_n]$, $n\ge 2$, be a homogeneous polynomial of degree $d$ defining 
a hypersurface $V(Q) \subset \PP^n$
having only isolated singularities. Let $F=F(Q)$ (respectively $M=M(Q)$) 
be the corresponding Milnor fiber (respectively complement), and let 
$p:F \to M$ be the canonical projection. 
The rational cohomology
$H^*(F,\Q)$ has a natural direct sum decomposition
$$H^*(F,\Q)=H^{*}(F,\Q)_{ 1}\oplus H^{*}(F, \Q)_{\ne 1}$$
where $H^{*}(F,\Q)_{ 1}=\ker (h^*-1)=p^*H^{*}(M,\Q)$ is the eigenspace corresponding to 
the eigenvalue $\lambda=1$ of the monodromy operator $h^*:H^{*}(F,\Q) \to H^{*}(F,\Q)$, 
and $H^{*}(F,\Q)_{\ne 1}=\ker((h^*)^{d-1}+...+1)$.
A different approach to the next result in the case of line arrangements can be found in \cite[(2.5.3)]{BDS}.

\begin{thm} \label{thm2}
With the above notation, the mixed Hodge structure on $H^{n-1}(F, \Q)$ is split, i.e.,
the subspaces $H^{n-1}(F,\Q)_{ 1}$ and $H^{n-1}(F, \Q)_{\ne 1}$
inherit pure Hodge structure from $H^{n-1}(F,\Q)$, such that $H^{n-1}(F,\Q)_{ 1}$ 
(respectively $H^{n-1}(F, \Q)_{\ne 1}$) has weight $n$ (respectively $n-1$). 
\end{thm}

\begin{proof}
To study the MHS on $H^*(M)$, one may use the following exact sequence of MHS, see \cite{PS}, p.138
\begin{equation} \label{MHS1}
...\to H^n(V(Q)) \to H^{n+1}_c(M) \to  H^{n+1}(\PP^n) \to ...
\end{equation}
Since $V(Q)$ has only isolated singularities, it follows that $ H^n(V(Q))$ 
is pure of weight $n$, see \cite{St}.
Moreover, the morphism $H^{n+1}(\PP^n) \to H^{n+1}(V(Q))$ is an isomorphism, hence $H^{n+1}_c(M)$
is pure of weight $n$. Using the duality between $H^*(M)$ and $H^*_c(M)$, 
see \cite{PS}, p.155, it follows that $H^{n-1}(M)$ is pure of weight $n$.

Let $\overline F$ be the projective hypersurface in $\PP^{n+1}$ defined by $\overline Q(x,t)=Q(x)+t^d=0$.
It is known that 
$\dim H^{n-1}(F)_{\ne 1}=\dim H^n(M(\overline Q))= \dim \coker (H^{n+1}(\PP^{n+1}) \to H^{n+1}(V(\overline Q))$,
see \cite{D1}, p.196 and p.206.

Next, the exact sequence of MHS
\begin{equation} \label{MHS2}
...\to H^n(V(Q)) \to H^{n+1}_c(F) \to  H^{n+1}(V(\overline Q)) \to H^{n+1}(V(Q))\to ...
\end{equation}
and the fact that  $ H^n(V(Q))$ (resp. $H^{n+1}(V(\overline Q))$) is pure of weight $n$ (resp. $n+1$), 
implies that the weight $n+1$ part in $H^{n+1}_c(F)$
has the same dimension as $\dim H^{n-1}(F)_{\ne 1}$. To see this, use the exact sequence
\begin{equation} \label{MHS3}
...\to  H^{n+1}(\PP^{n+1})   \to H^{n+1}(V(\overline Q) ) \to H^{n+2}_c(M( \overline Q ) ) \to  H^{n+2}(\PP^{n+1}) \to ...
\end{equation}
By duality,  $\dim W_{n-1}H^{n-1}(F)=\dim H^{n-1}(F)_{\ne 1}$.

Now, the equality $H^{n-1}(F,\Q)_{ 1}=p^*H^{n-1}(M,\Q)$ and our proof above show that
$H^{n-1}(F,\Q)_{ 1}$ is a pure HS of weight $n$. 

Consider now the subspace $E=W_{n-1}H^{n-1}(F,\Q)$.
Then $E$ is $h^*$ invariant, $E \cap H^{n-1}(F,\Q)_{ 1}=0$ and $\dim E= \dim H^{n-1}(F,\Q)_{\ne 1}$.
These three properties imply that $E= H^{n-1}(F)_{\ne 1}$. Since the proof above implies also that 
$W_{n-2}H^{n-1}(F,\Q)= W_{n+2}H^{n+1}_c(F,\Q)= 0$, all the claims are proven.
\end{proof}

As an application, we offer a converse to Corollary \ref{cor:monores}, for $q=1$.

\begin{cor}
\label{cor:rdetect}
Let $p:F\to M$ be the Galois cover associated to an arbitrary line arrangement.
The algebraic monodromy action on $H^1(F)$ is trivial if and only if 
$p^*(\RR_1(M))=\RR_1(F)$.
\end{cor}

\begin{proof}
We will show that $\RR_1(F) \not\subseteq p^*(\RR_1(M))$, if $h^* \ne \id$ on $H^1(F)$.
Indeed, we know from Theorem \ref{thm2} that 
$W_1H^1(F)= \bigoplus_{\lambda \ne 1} H^1(F)_{\lambda} \ne 0$.
Let us first assume that $ H^1(F)_{\lambda} \ne 0$, for some $\lambda \ne \pm 1$.
Pick linearly independent eigenvectors of $h^*$, $\al \in H^1(F)_{\lambda}$ and
$\be \in H^1(F)_{\lambda^{-1}}$. We will conclude in this case by showing that
$\al \cup \be =0$. Clearly, $\al \cup \be \in W_2H^2(F)\cap p^*(H^2(M))= p^*(W_2H^2(M))$.
Since the MHS on $H^2(M)$ is pure of type $(2,2)$, by \cite{Sh}, $W_2H^2(M)=0$, and
we are done. In the remaining case, $W_1H^1(F)= H^1(F)_{-1} \ne 0$. Since this
space has even dimension, we may find two linearly independent eigenvectors of $h^*$, 
$\al ,\be \in H^1(F)_{-1}$. Then $\al \cup \be \in  p^*(W_2H^2(M))=0$, as before.
\end{proof}

We close this section with an addendum to Theorem \ref{thm:monoch}, concerning
translated components. Let $M=M(Q)$ be a curve complement in $\PP^2$,
with associated cyclic Galois cover $p:F\to M$. Assume $h^*=\id$ on $H^1(F)$. 

We need to briefly recall the general theory, following \cite{A} and \cite{D3}. 
Let $X$ be a smooth, irreducible, quasi-projective variety. The positive-dimensional
components, $W'$, of $\V_1(X)$, are intimately related to {\em admissible maps},
that is, surjective regular maps onto smooth curves, $f:X\to S$, having connected generic fiber.
Each $W'$ is a (possibly translated) subtorus of the form $\rho f^*(\T(S))$, 
where $f$ is admissible, $\chi (S)\le 0$, and $\rho \in \T(X)$ has finite order.
The number of components with the same {\em direction}, $W:= f^*(\T(S))$, 
denoted by $n(W)$, is computable in terms of the {\em multiplicities} of the fibers of $f$.
For $c\in S$, we denote the corresponding multiplicity by $m_c(f)$.

A key result due to Arapura, Proposition V.1.7 from \cite{A}, guarantees that
$f^*(\T(S))$ is a component of $\V_1(X)$, whenever $f$ is admissible and $\chi (S)<0$.
We will need the following converse implication.

\begin{prop}
\label{prop:conv}
Let $f:X\to S$ be a surjective regular map. If $\chi (S)<0$ and $f^*(\T(S))$ is 
an irreducible component of $\V_1(X)$, then $f$ is admissible.
\end{prop}

\begin{proof}
Take a Stein factorisation, $f=q\circ f'$, where $f': X\to S'$ is admissible, and
$q:S'  \to S$ is finite. Now, we will make an estimate for Euler numbers, valid for
an arbitrary finite map $q$ between smooth, irreducible curves.

Let $e$ be the degree of $q$. Let $S =\cup_{k=1,e}S_k$ be the canonical partition:
each $S_k$  consists of those $c\in S$ having exactly $k$ $q$-preimages. For $k<e$,
$S_k$ is a finite set, of cardinality $s_k\ge 0$. Then $S'=\cup_{k=1,e} q^{-1} (S_k)$,
and the restriction $q: q^{-1}(S_e)\to S_e$ is a non-ramified $e$-fold cover. Using
the additivity of  Euler numbers with respect to constructible partitions (see \cite{D2}),
we find that
\begin{equation} \label{eq:chi1}
\chi (S')-\chi (S)= (e-1)\chi (S_e)+ \sum_{1\le k<e} (k-1)s_k\, ,
\end{equation}
which leads to
\begin{equation} \label{eq:chi2}
\chi (S')-\chi (S)= (e-1)\chi (S)+ \sum_{1\le k<e} (k-e)s_k\, .
\end{equation}

We infer from \eqref{eq:chi2} that $\chi (S')<0$, since $\chi (S)<0$. Therefore,
$f'^*(\T(S'))$ is an irreducible component of $\V_1(X)$, containing $f^*(\T(S))$,
which implies that $f'^*(\T(S'))=f^*(\T(S))$. In particular $b_1(S')=b_1(S)$,
hence $\chi (S')=\chi(S)$. Since $(k-e)s_k \le 0$ for any $k$, and $\chi(S)<0$, 
we deduce from \eqref{eq:chi2} that $e=1$. Therefore, the generic fiber of $f$ 
is connected, as asserted, since $f'$ has this property.
\end{proof}

Going back to the case when $X=M(Q)$ is a curve complement in $\PP^2$, and
$h^*=\id$ on $H^1(F)$, we recall from Theorem \ref{thm:monoch} the correspondence
$W\mapsto p^*(W)$. This gives a dimension-preserving bijection between non-translated
subtori in $\T(M)$ and $\T(F)$, that identifies the non-translated irreducible
components in $\V_1(M)$ and $\V_1(F)$. In geometric terms, this may be rephrased as follows.

\begin{cor}
\label{cor:geombij}
Let $p:F\to M$ be the Galois cover associated to a plane projective curve.
Assume the algebraic monodromy action on $H^1(F)$ is trivial. If $f:M\to S$ is admissible
and $\chi(S)<0$, then $g:= f\circ p: F\to S$ is  admissible, too.
\end{cor}

Note also that any admissible map, $f:M\to S$ or $g:F\to S$, must be a rational pencil
(compare to Corollary \ref{cor:subtle}). This is due to the fact that $W_1H^1(M)=0$,
and consequently $W_1H^1(F)=0$; see \cite[Proposition 7.2]{DPS08}.

\begin{prop} \label{prop1}
Under the assumptions from Corollary \ref{cor:geombij}, the following hold.

(i) The bijection $W\mapsto p^*(W)$ identifies, for $k>1$, the $k$-dimensional directions
of irreducible components in $\V_1(M)$ and $\V_1(F)$. Moreover, $n(W)=n(p^*(W))$.

(ii) If $W$ is a $1$-dimensional direction
of irreducible component in $\V_1(M)$, then so is $p^*(W)$ in $\V_1(F)$, and
$n(W)\le n(p^*(W))$.
\end{prop}

\begin{proof}
$(i)$ For $k>1$, it follows from \cite{D3} that (for both $M$ and $F$) the
$k$-dimensional directions from our statement coincide with the non-translated
$k$-dimensional components of $\V_1$. Moreover, they are all of the form $f^*(\T(S))$,
where $f$ is an admissible map onto a rational curve with $\chi(S)=1-k$, and 
$n (f^*(\T(S)))= \prod_c m_c(f)$, where the product is taken over those $c\in S$ with 
$m_c(f)>1$. Our claims follow then from Theorem \ref{thm:monoch}, Corollary \ref{cor:geombij},
and the remark that $m_c(f\circ p)=m_c(f)$, for all $c\in S$, since $p$ is a submersion.

$(ii)$ For $k=1$, we know from \cite{D3} that the $1$-dimensional directions 
(for both $M$ and $F$) coincide with the subtori of the form $f^*(\T(\C^*))$,
where $f$ is an admissible map onto $\C^*$, having at least one multiple fiber.
In this case, $1+n (f^*(\T(\C^*)))= \prod_c m_c(f)$. If $f:M\to \C^*$ satisfies the above
conditions, let us consider the map $g:=f\circ p$.

We may factor it as in the proof of Proposition \ref{prop:conv}, $g=q\circ g'$,
with $g':F\to S'$ admissible and $q:S'\to \C^*$ finite. We also know from
\cite[Corollary 3.21]{D-p} that $f^*(H^1(\C^*, \C))$ is a maximal isotropic subspace in
$H^1(M)$, which implies the same property, for $g^*(H^1(\C^*, \C))$ in $H^1(F)$.
The inclusion $g^*(H^1(\C^*, \C))\subseteq g'^*(H^1(S', \C))$ must then be an equality,
which forces $S'$ to be $\C^*$, and $q$ to induce a cohomology isomorphism. In this case,
\eqref{eq:chi2} becomes 
\[
\sum_{k<e} (k-e)s_k =0\, ,
\]
therefore $q$ is unramified. This implies that 
$m_c(f)= m_c(q\circ g')= {\rm g.c.d.}~ \{ m_{c'} (g') \mid c'\in q^{-1}(c) \}$, for all $c\in \C^*$.

We thus see that each multiple fiber $f^{-1}(c)$ gives rise to $e$ multiple fibers
of $g'$. In particular, $g'^*(\T(\C^*))= g'^* \circ q^*(\T(\C^*))= p^* \circ f^*(\T(\C^*))$
is a $1$-dimensional direction in $\V_1(F)$. Clearly, 
$ \prod_c m_c(f)^e \le  \prod_{c'} m_{c'}(g')$, which verifies our last claim.
\end{proof}


\begin{thebibliography}{00}

\bibitem{A} D. Arapura, 
{\em Geometry of cohomology support loci for local systems.} I, 
J. Algebraic Geom. \textbf{6} (1997), 563--597.

\bibitem{Beau} A.~Beauville, 
{\em Annulation du $H\sp 1$ pour les fibr\' es en droites plats}, in: 
{\em Complex algebraic varieties} (Bayreuth, 1990), 
Lect. Notes in Math. \textbf{1507}, Springer, Berlin, 1992, pp. 1--15. 

\bibitem{BP} B.~Berceanu, S.~Papadima,
{\em Universal representations of braid and braid-permutation 
groups}, preprint {\tt arXiv:0708.0634}, to appear in 
J. Knot Theory Ramifications.

\bibitem{BG} A.~K.~Bousfield, V.~K.~A.~M.~Gugenheim,
{\em On} PL {\em De Rham theory and rational homotopy type},
Memoirs Amer. Math. Soc., vol.~179, Amer. Math. Soc., 
Providence, RI, 1976.

\bibitem{Br} E.~Brieskorn,
{\em Sur les groupes de tresses}, in:
S\'{e}minaire Bourbaki, 1971/72, Lect. Notes in Math. 
\textbf{317}, Springer-Verlag, 1973, pp.~21--44. 

\bibitem{B}  K.~S.~Brown,
{\em Cohomology of groups}, Grad. Texts in Math., 
vol.~87, Springer-Verlag, New York-Berlin, 1982.

\bibitem{BDS} N. Budur, A. Dimca, M. Saito, {\em First Milnor cohomology of hyperplane arrangements},  preprint {\tt arXiv:0905.1284}.

\bibitem{DP} A.~D.~R.~Choudary, A.~Dimca, S.~Papadima,
{\em Some analogs of Zariski's theorem on nodal line arrangements},
Algebraic and Geometric Topology \textbf{5} (2005), 691--711.

\bibitem{CS1} D.~C.~Cohen, A.~I.~Suciu,
{\em On Milnor fibrations of arrangements}, 
J. London Math. Soc. \textbf{51} (1995), no.~2, 105--119.

\bibitem{DGMS}  P.~Deligne, P.~Griffiths, J.~Morgan, D.~Sullivan,
{\em Real homotopy theory of {K}\"{a}hler manifolds},
Invent. Math. \textbf{29} (1975), no.~3, 245--274.

\bibitem{D1} A. Dimca, {\em Singularities and Topology of Hypersurfaces}, 
Universitext, Springer-Verlag, 1992.

\bibitem{D2} A. Dimca,
{\em Sheaves in Topology},  Universitext, Springer-Verlag, 2004.

\bibitem{D3} A. Dimca,
{\em Characteristic varieties and constructible sheaves}, 
Rend. Lincei Mat. Appl. \textbf{18} (2007),
365--389.

\bibitem{D4} A. Dimca,
{\em On the isotropic subspace theorems}, 
Bull. Math. Soc. Sci. Math. Roumanie, \textbf{51} (2008), 
no.~4, 307--324. 

\bibitem{D-p} A. Dimca,
{\em Pencils of plane curves and characteristic varieties},
preprint {\tt math.AG/0606442}.

\bibitem{D5} A. Dimca,
{\em On admissible rank one local systems}, J. of Algebra, \textbf{321} (2009), 3145--3157.  


\bibitem{DPSqk}  A.~Dimca, S.~Papadima, A.~Suciu, 
{\em Quasi-K\" ahler groups, $3$-manifold groups, and formality}, preprint
{\tt arXiv:0810.2158}.

\bibitem{DPS08} A.~Dimca, S.~Papadima, A.~Suciu,
{\em Topology and geometry of cohomology jump loci}, 
preprint {\tt arXiv:0902.1250}, to appear in Duke Math. J.



\bibitem{FY} M. Falk, S. Yuzvinsky,
{\em Multinets, resonance varieties, and pencils of plane curves}, 
Compositio Math. \textbf{143} (2007), no.~4, 1069--1088.

\bibitem{GL} M.~Green, R.~Lazarsfeld, 
{\em Higher obstructions to deforming cohomology groups of 
line bundles}, J. Amer. Math. Soc. \textbf{4} (1991), no.~1, 87--103.


\bibitem{KM} M.~Kapovich, J.~Millson, 
{\em On representation varieties of {A}rtin groups, projective 
arrangements and the fundamental groups of smooth complex 
algebraic varieties}, Inst. Hautes \'{E}tudes Sci. Publ. Math.
\textbf{88} (1998), 5--95. 


\bibitem{L1} A.~Libgober,
{\em Eigenvalues for the monodromy of the Milnor fibers of arrangements},
in: {\em Trends in singularities}, 
Trends Math., Birkh\" auser, Basel, 2002, pp. 141--150.

\bibitem{Mac}  A.~Macinic,
{\em Cohomology rings and formality properties of nilpotent 
groups}, preprint {\tt arXiv:0801.4847}.

\bibitem{MP}  A.~Macinic, S.~Papadima,  
{\em On the monodromy action on {M}ilnor fibers 
of graphic arrangements}, Topology Appl., 
\textbf{156} (2009), 761--774. 

\bibitem{M}  J.~W.~Morgan,
{\em The algebraic topology of smooth algebraic varieties},
Inst. Hautes \'{E}tudes Sci. Publ. Math. \textbf{48} (1978), 137--204.

\bibitem{OS} P.~Orlik, L.~Solomon,
{\em Combinatorics and topology of complements of
hyperplanes}, Invent. Math. \textbf{56} (1980), 167--189.

\bibitem{PS-t} S.~Papadima, A.~I.~Suciu, 
{\em Toric complexes and {A}rtin kernels}, Advances in Math. 
\textbf{220} (2009), no.~2, 441--477.

\bibitem{PS-m} S.~Papadima, A.~I.~Suciu, 
{\em Algebraic monodromy and obstructions to formality}, preprint
{\tt arXiv:0901.0105}, to appear in Forum Math.

\bibitem{PSs} S.~Papadima, A.~I.~Suciu, {\em Geometric and algebraic aspects of 1-formality},
preprint {\tt arXiv:0903.2307}, to appear in Bull. Math. Soc. Sci. Math. Roumanie.


\bibitem{PS} C. Peters, J. Steenbrink, {\em Mixed Hodge Structures}, 
Ergeb. der Math. und ihrer Grenz. 3. Folge 52,
Springer, 2008.

\bibitem{Q} D.~Quillen,
{\em Rational homotopy theory}, Ann. of Math.
\textbf{90} (1969), no.~2, 205--295.  

\bibitem{Sh} B.~Z.~Shapiro,
{\em The mixed Hodge structure of the complement to an arbitrary arrangement of affine
complex hyperplanes is pure}, Proc. Amer. Math. Soc. \textbf{117} (1993), no.~4, 931--933.

\bibitem{St} J. Steenbrink, {\em Mixed Hodge structures associated with isolated singularities}, 
in: {\em Singularities, Part 2} (Arcata, 1981), Proc. Symp. Pure Math. \textbf{40}, 
Amer. Math. Soc., 1983, pp. 513-536.

\bibitem{S2} A. Suciu,
{\em Fundamental groups of line arrangements: enumerative aspects}, in: 
{\em Advances in algebraic geometry motivated by physics} (Lowell, MA, 2000), 
Contemp. Math. \textbf{276}, Amer. Math. Soc., Providence, RI, 2001, pp. 43--79. 

\bibitem{Sul}  D.~Sullivan,
{\em Infinitesimal computations in topology},
Inst. Hautes \'{E}tudes Sci. Publ. Math.
\textbf{47} (1977), 269--331.

\bibitem{Y} S.~Yuzvinsky,
{\em A new bound on the number of special fibers in a pencil of curves}, Proc.\ Amer.\ Math.\ Soc.\ \textbf{137} (2009), 1641--1648.


\end{thebibliography}
\end{document}